\documentclass{amsart}
\usepackage{hyperref}
\usepackage{algpseudocode}
\usepackage{graphicx}

\hypersetup{
    pdftitle={An Improved Algorithm for hypot(a,b)},
    pdfauthor={Carlos F. Borges},
    pdfkeywords={hypot(), floating point, fused multiply-add, IEEE 754}
    }

\newtheorem{lemma}{Lemma}
\newtheorem{algorithm}{Algorithm}
\numberwithin{equation}{section}

\def\F{\mathbb{F}}
\def\Re{\mathbb{R}}

\begin{document}

\title{An Improved Algorithm for {\tt hypot(a,b)}} 

\author{Carlos F. Borges}
\address{Department of Applied Mathematics\\Naval Postgraduate School\\Monterey CA 93943}
\email{borges@nps.edu}

\begin{abstract}
We develop a fast and accurate algorithm for evaluating $\sqrt{a^2+b^2}$ for two floating point numbers $a$ and $b$. Library functions that perform this computation are generally named {\tt hypot(a,b)}. We will compare four approaches that we will develop in this paper to the current resident library function that is delivered with Julia 1.1 and to the code that has been distributed with the C math library for decades. We will demonstrate the performance of our algorithms by simulation.
\end{abstract}


\date{April 30, 2019.}


\maketitle

Given two floating point numbers $a$ and $b$ we wish to compute $\sqrt{a^2+b^2}$ quickly and accurately. Library functions that perform this computation are generally named {\tt hypot(a,b)}. And although the IEEE 754 standard suggests (but does not require) that computer languages provide a version that gives a correctly rounded answer, we are not aware of any current language in wide use that does so. 

We note that the problem is trivial if either operand is zero and further that the magnitudes of the operands are important but the signs are not. Therefore, we shall confine our in depth analysis to the case where $a \geq b > 0$. Enforcing this condition will be an initialization step when we develop an actual code.  This paper will be restricted to the case where all floating point calculations are done in IEEE 754 compliant arithmetic using round-to-nearest rounding, although many of the results can be extended to other formats under proper conditions.

\section{Existing approaches}

The need to compute $\sqrt{a^2+b^2}$ is common in numerical computation and so there are a few existing standard codes for doing so. A common textbook approach is the widely known {\em stable} algorithm for this computation that is used in Julia 1.1 which we shall call {\tt Julia1.1\_hypot(a,b)}. It essentially amounts to 

\vspace{.2in}
\begin{algorithm} {\tt Julia1.1\_hypot()} 
\hrule
\begin{algorithmic}
\If{{\tt a == 0}}
	\State {\tt h = 0}
\Else 
	\State {\tt r = b/a}
	\State {\tt h = a*sqrt(1+r*r)}
\EndIf

\end{algorithmic}
\hrule
\end{algorithm}
\vspace{.2in}

This approach avoids unecessary overflow/underflow that might occur in an interim calculation when very large/small inputs are squared. In the non-trivial case it accomplishes this by rescaling so that the quantity $r$ that must be squared is $1 \geq r \geq 0$ and this cannot lead to an overflow/underflow error prior to the calculation of the square root. Overflow is only possible if the true value of $\sqrt{a^2+b^2}$ is beyond the range of the floating point format. 

The downside of this approach is that it rescales even when rescaling is unecessary and that rescaling generates additional rounding error. We will see in our experiements that this approach sometimes give answers that are as far as two ulps from the correctly rounded answer.

The second existing approach is the code that has been delivered with the standard C math library since at least 1991. The code appears in appendix B and we shall refer to it as {\tt clib\_hypot()}. It deals with widely varying operands and prevents intermediate overflow/underflow using methods similar to those we will use in ours.\footnote{There is a clear flaw in the threshold used for widely varying operands and it is much broader than necessary.} It works by using a variety of tricks to artificially extend the precision of the computed $a^2+b^2$ so that they get a correctly rounded value (or very nearly so). This is then passed to the {\tt sqrt()} function. What we will observe in the testing is that all this additional work can be effectively superseded by a single fused multiply-add.

\section{Mathematical Preliminaries}

We begin with a few definitions. We shall denote by $\F \subset \Re$ the set of all strictly positive floating point numbers in the normalized range of the current format. We define $fl(x):\Re \rightarrow \F$ to be a function such that $fl(x)$ is the element of $\F$ that is closest to $x$ provided $x$ lies within the normalized range of the current format. We define {\em machine epsilon}, which we will denote by $\epsilon$, to be the difference between $1$ and the closest larger element of $\F$. For example, in IEEE 754 double precision $\epsilon = 2^{-52}$. We define $F_{max}$ and $F_{min}$ to be, respectively, the largest and smallest elements of $\F$.

We will also need the following three lemmas. The first two both relate to floating point computations.

\begin{lemma}
\label{fpl1}
Let $x \in \Re$ lie within the normalized range of the current format. Then
$$
fl(x) = x (1+\delta)
$$
for some $\delta$ satisfying $|\delta|<\frac{\epsilon}{2} $
\end{lemma}

This lemma is very well known and can be found in \cite{Overton}. A lesser known but related lemma is:\footnote{To see why it's true simply consider the operation {\tt 1.0*(1.0 - eps/3.9)}.}

\begin{lemma}
\label{fpl2}
Let $x \in \F$ be any normalized floating point number. Then
$$
x = fl(x(1+\delta))
$$
for any $\delta$ satisfying $|\delta| < \frac{\epsilon}{4}$.
\end{lemma}

And finally, this analytical lemma will be useful in what follows and is easily established by examining the Maclaurin series expansion of $\sqrt{1+x}$:

\begin{lemma}
If $x > 0$ then $\sqrt{1+x} < 1 + \frac{x}{2}$
\label{l1}
\end{lemma}

\section{Operands with widely differing magnitudes}

To see what happens when $a > b$ have widely differing magnitudes observe that
\begin{eqnarray*}
\sqrt{a^2+b^2} & = & a \sqrt{1 + (b/a)^2} \\
                       & < & a \left(1 + \frac{(b/a)^2}{2}\right)
\end{eqnarray*}
where the inequality follows by replacing $x$ with $(b/a)^2$ in lemma \ref{l1}. Now, if 
$$
\frac{(b/a)^2}{2} \leq \frac{\epsilon}{4}
$$
then lemma \ref{fpl2} guarantees that $a = fl(\sqrt{a^2+b^2})$. We can rewrite the inequality and accept $a$ as the correctly rounded answer whenever $b \leq a\sqrt{\epsilon/2}$. Note that this form will also yield the exact answer whenever $b=0$. If we test for this case first in our code then we will not have to separately check for zero operands in the initialization phase.

\section{Operands without widely differing magnitudes}
\label{notwide}

When the operands do not have widely differing magnitudes then our approach will be to compute $a^2+b^2$ with some care and use the built-in {\tt sqrt()} function. There could be an unecessary floating point exception if the calculation of $a^2+b^2$ leads to an overflow/underflow/denormalization and we need to avoid that where possible. Clearly, there can be no overflow if $a \leq \sqrt{F_{max}/2}$ and there can be no underflow/denormalization provided that $b \geq \sqrt{F_{min}}$. If one of these is violated\footnote{At this stage of the process it is not possible to violate both because such a pair of operands would have widely differing magnitudes.} we can avoid the exception by rescaling the operands, doing the calculation, and then scaling back. When such rescaling is necessary there are two options. Some codes use the standard library functions {\tt frexp()} and {\tt ldexp()} to directly change the exponents and hence effect the rescaling without incurring rounding errors. The only downside of this is that it can be a bit more costly and it is not as direct. A better approach is to rescale by multiplying by a power of the base. This has the same effect as directly changing the exponents (i.e. there is no rounding error) but it leverages the hardware floating point operations which are generally much faster. For this to work the rescaling constant needs to be an appropriate integer power of the base. The precise range of appropriate values is dependent on the particular floating point format but one relatively format independent method for finding an appropriate rescaling value is to use {\tt eps(sqrt(floatmin(T)))} where {\tt T} is the floating point type. This will give a rescaling whose value is smaller than $\sqrt{F_{min}}$ and whose reciprocal is larger than $\sqrt{F_{max}/2}$. This can be used to rescale and then scale back the operands with no rounding error. Moreover, the cumbersome calculation of this constant should be taken care of at compile time so there is no additional cost.

A little rounding error analysis will be very useful at this point. Let $z \in \F$ be the result of our floating point calculation of $a^2+b^2$. Because this calculation will be done in floating point there must be some $\delta_1$, whose value we do not yet know, such that
$$
z = (a^2+b^2)(1+\delta_1).
$$

Now let $h = \mbox{{\tt sqrt}}(z)$. Since we are assuming that the {\tt sqrt()} function is correctly rounded lemma \ref{fpl1} implies that $h = \sqrt{z} (1+\delta_2)$ for some $|\delta_2| < \frac{\epsilon}{2}$. Putting this all together and invoking lemma \ref{l1} leads to

\begin{eqnarray*}
h & = & \sqrt{(a^2+b^2)(1+\delta_1)} (1+\delta_2) \\
   & < & \sqrt{a^2+b^2} (1+\frac{\delta_1}{2}) (1+\delta_2) \\
   & < & \sqrt{a^2+b^2} (1+\frac{\delta_1+2\delta_2}{2} +\frac{\delta_1\delta_2}{2})
\end{eqnarray*}

The first thing we should observe is that if we were able to compute a correctly rounded $z$ then $|\delta_1| < \frac{\epsilon}{2}$. Plugging in above we see that the relative error in the computed $h$ could be as large as
$$
\frac{3\epsilon}{4} +\frac{\epsilon^2}{8}
$$
and hence the computed $h$, after rounding, could be as much as one $ulp(h)$ different from the correctly rounded true value.

In a similar fashion if $|\delta_1| < \frac{3\epsilon}{2}$ then the relative error in the computed $h$ could be as large as
$$
\frac{5\epsilon}{4} +\frac{3 \epsilon^2}{8}
$$
and hence the computed $h$, after rounding, could not be more than one $ulp(h)$ different from the correctly rounded true value. Careful rounding error analysis of the computation $z = a*a+b*b$ shows that the computed value satisfies
$$
|z-(a^2+b^2)| \leq \epsilon z
$$
and hence our computed answer must be within one ulp of the correctly rounded one. An excellent treatment can be found in section 5.3 of \cite{Jrod}.

\subsection{Naive Approaches}

At this point we propose two naive algorithms for the computation in the event that the operands are not widely separated. First is the direct calculation
\vspace{.2in}
\begin{algorithm} Naive (Unfused)

\hrule
\begin{algorithmic}
	\State {\tt h = sqrt(a*a+b*b)}
\end{algorithmic}
\hrule
\end{algorithm}
\vspace{.2in}

And the second uses the fused multiply-add operation which allows us to get a more accurate calculation of $a^2+b^2$.

\vspace{.2in}
\begin{algorithm} Naive (Fused)

\hrule
\begin{algorithmic}
	\State {\tt h = sqrt(fma(a,a,b*b))}
\end{algorithmic}
\hrule
\end{algorithm}
\vspace{.2in}

\subsection{Differential Correction}

As we noted at the end of section \ref{notwide} the correctly rounded square root of the correctly rounded $a^2+b^2$ can still be off by as much as one ulp. This hints at the possibility that working harder to compute $a^2+b^2$ more accurately may not be the best path to a better answer. That leads us to this final approach which we shall see has a great deal of merit.
 
A subtle and highly underappreciated phenomenon of floating point computation is that formally invertible functions are often not one-to-one in correctly rounded floating point. The square root function is a perfect example. If we let $S = \left\{ x \in \F | 1 \leq x < 4\right\}$ then it is clear that {\tt sqrt()} maps the finite set $S$ onto a proper subset of itself and hence is not one-to-one by the pigeonhole principle.

Why is this is important? Because it tells us that we should not be looking for a correction $x$ that 
$$
a^2+b^2+x = fl(a^2+b^2)
$$
but rather we should be looking for a correction $x$ such that if $h = sqrt(a*a+b*b)$ is the computed hypoteneuse we have
$$
a^2+b^2+x = h^2
$$
If we had such an $x$ then we could use the Maclaurin series expansion
$$
\sqrt{h^2-x} = h - \frac{x}{2h} + O(x^2)
$$
 truncated to first order to correct our computed hypoteneuse value.\footnote{We note that our testing indicates that there is no point in taking the Maclaurin series correction past the first order term.}

Setting $h^2-x = a^2+b^2$ allows us to find that $x = h^2-a^2-b^2$ and although we cannot generally compute $x$ exactly, we do want to do so very carefully. If we let $\delta = h-a$ then the reader can verify that
\begin{eqnarray*}
x & = & h^2 - a^2 - b^2\\
& = & 2\delta(a-b) + (2\delta-b)b + \delta^2
\end{eqnarray*}
Note that because $h$ and $a$ are nearly equal numbers there is no rounding error in computing $\delta = h-a$. The precise arrangement and evaluation order of this formula is critical to get the best performance over the full range of values. 

The form above gives good performance over the full range of operands but we can do better if we switch approaches based on the sizes of the arguments. In particular, if $h \leq 2b$ then we will set $\delta = h-b$ and use
$$
a(2\delta - a)+(\delta-2(a-b))\delta
$$
otherwise, we will set  $\delta = h-a$ and use
$$ 
2\delta(a-2b) + (4\delta-b)b + \delta^2
$$
in both cases $\delta$ is computed without rounding error. It is worth noting that there are many other ways to compute $x$ and one can break this problem into lots of tiny pieces and do better but we feel that this approach is the best tradeoff between accuracy and efficient code. We can now use this to apply a correction to the computed value of $h$ which leads to our third algorithm which uses only the standard floating point operations.

\vspace{.2in}
\begin{algorithm} Corrected (Unfused)

\hrule
\begin{algorithmic}
\State {\tt h = sqrt(a*a+b*b)}
\If{{\tt h <= 2*ay}}
	\State {\tt delta = h-ay}
	\State{\tt x = ax*(2*delta - ax)+(delta-2*(ax-ay))*delta}
\Else 
	\State {\tt delta = h-ax}
	\State {\tt x = 2*delta*(ax-2*ay) + (4*delta-ay)*ay + delta*delta}
\EndIf
\State {\tt h = h - x/(2*h)}
\end{algorithmic}
\hrule
\end{algorithm}
\vspace{.2in}

Finally, if it happens that we are on an architecture that supports the fused multiply-add then it is possible to compute the correction far more accurately. This leads to our final code:

\vspace{.2in}
\begin{algorithm} Corrected (Fused)

\hrule
\begin{algorithmic}
\State {\tt h = sqrt(fma(a,a,b*b))}
\State {\tt h\_sq = h*h}
\State {\tt a\_sq = a*a}
\State {\tt x = fma(-b,b,h\_sq-a\_sq) + fma(h,h,-h\_sq) - fma(a,a,-a\_sq)}
\State {\tt h = h - x/(2*h)}
\end{algorithmic}
\hrule
\end{algorithm}
\vspace{.2in}

\section{Testing}

We now test our proposed algorithms against the two existing approaches. As a baseline for testing purposes we will use the big float format in Julia to do an arbitrary precision calculation of the hypoteneuse and then cast the result to a double precision floating point number to get a {\em correctly rounded value} that we shall call $\bar{h}$. We note that in some exceedingly rare situations it may not be true that this yields the correctly rounded value (due to the effects of double rounding) but it is certainly good enough for our purposes.

Our first test will consist of creating $10^9$ random pairs of double precision floating point numbers where both are distributed according to a standard normal distribution. We will run all of our algorithms on each pair as well as computing $\bar{h}$. We summarize the percentage of times each algorithm differed from $\bar{h}$ by exactly one ulp and exactly two ulps in table \ref{tab:table1}.

\begin{table}[h]
  \begin{center}
    \caption{Testing done with $a,b \sim \mathcal{N}(0,1)$}
    \label{tab:table1}
    \begin{tabular}{l|c|c}
      \textbf{Method} & \textbf{One ulp errors (\%)} &  \textbf{Two ulp errors (\%)}\\
      \hline
      {\tt Julia1.1\_hypot()} & 35.08 & 0.16\\
      {\tt clib\_hypot()} & 12.91 & 0\\
      Naive (Unfused) & 16.70 & 0\\
      Naive (Fused) & 13.02 & 0\\
      Corrected (Unfused) & 0.54 & 0\\
      Corrected (Fused) & 0 & 0\\
    \end{tabular}
  \end{center}
\end{table}

We immediately observe several things. First of all, the {\tt Julia1.1\_hypot()} code which is a very widely known approach to the problem performs rather horribly. The willingness to add rounding error to escape intermediate overflow/underflow problems destroys accuracy. Not only does it have the worst performance of all the approaches it is the only one that can produce two ulp errors (although this is admittedly rare). Second, the {\tt clib\_hypot()} and the Naive (Fused) approaches yield very comparable results. Both are basically trying to improve accuracy by computing $a^2+b^2$ with additional precision and so it is not surprising that they behave similarly even though the get the additional precision in different ways. Unfortunately, both suffer from the effects of iterated rounding. Finally, and most importantly, we see very clearly the superiority of the corrected approach (whether fused or unfused). When leveraging the fused multiply-add we get correctly rounded results every time. And the unfused version runs almost as well.

\begin{figure}[h] 
\begin{center}
\includegraphics[width=6in]{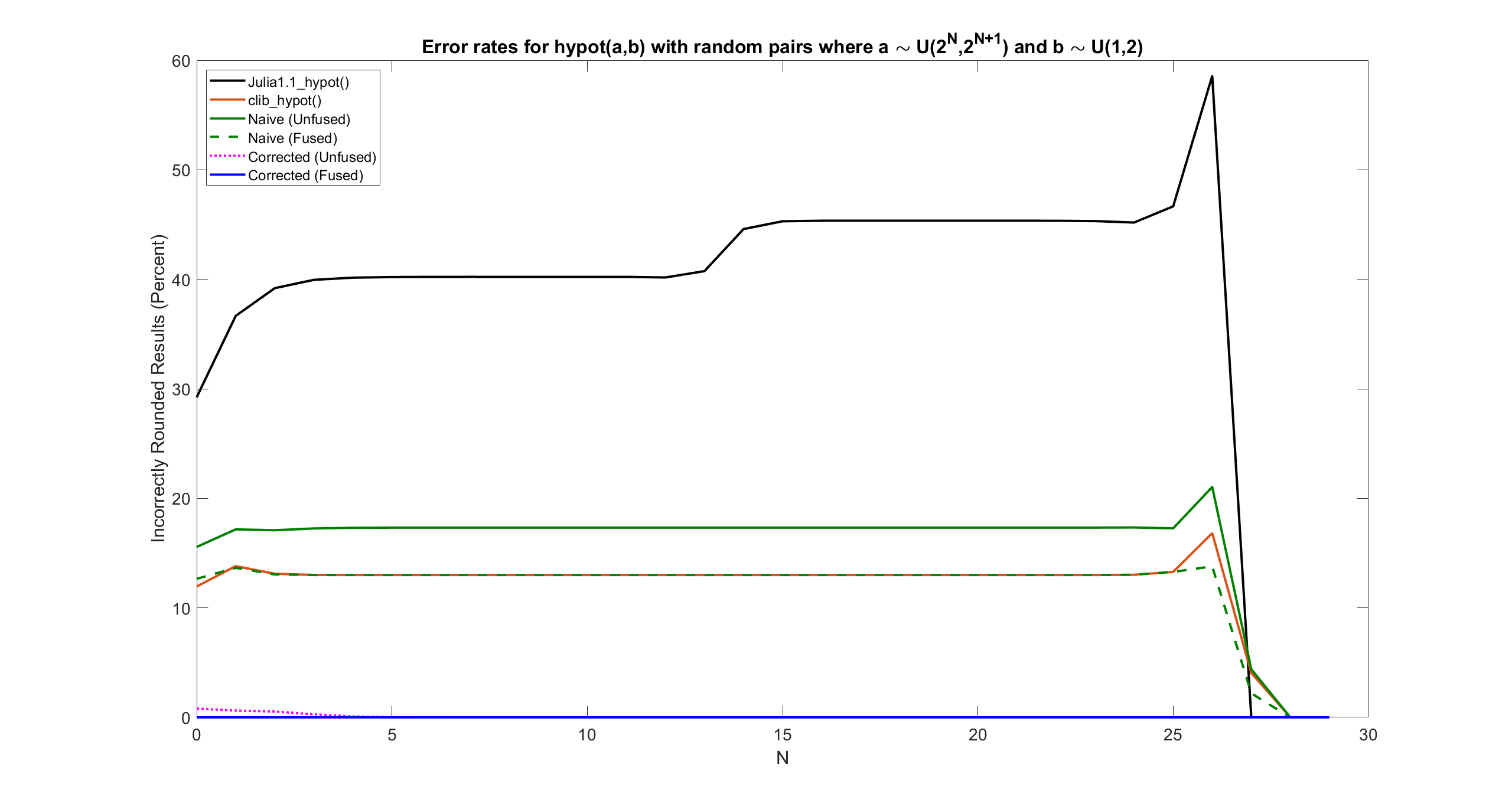}
\caption{Comparison of the six algorithms over a range of relative scales.}
\label{bigtestplot}
\end{center}
\end{figure}

A further experiment will more completely display the differences in these algoritms. This time we track their accuracy over a range of relative scales of the operands. For this experiment we will use uniformly distributed operands. In particular, we will let $a \sim U(2^N,2^{N+1})$ and $b \sim U(1,2)$. In effect, our test pairs will be floating point numbers whose significands are uniformly distributed, but whose exponents differ by $N$. We will run $10^9$ random pairs for each value $N=0,1,...,29$ and will compute the percentage of incorrectly rounded results for each of the six approaches. The results appear graphically in figure \ref{bigtestplot} as well as numerically in table \ref{tab:table2}. One can see the stunning difference in performance between the approaches. Note that in table \ref{tab:table2} we observe that we are getting true correctly rounded results more than $99\%$ of the time for all relative scales for the unfused version of the corrected approach.

\begin{table}[p]
  \begin{center}
    \caption{Percentage of incorrectly rounded results by relative scale for  $a \sim U(2^N,2^{N+1})$ and $b \sim U(1,2)$}
    \label{tab:table2}
    \begin{tabular}{l|l|l|l|l|l|l}
        &  &   & \multicolumn{2}{c}{Naive} & \multicolumn{2}{c}{Corrected} \\
     N & {\tt Julia1.1\_hypot()} &  {\tt clib\_hypot()} & Unfused & Fused & Unfused & Fused \\
      \hline
0 & 29.2282102 & 11.9488456 & 15.5629732 & 12.6486134 & 0.8005489 & 0\\
1 & 36.6616214 & 13.8082644 & 17.1661734 & 13.6513319 & 0.6212053 & 0\\
2 & 39.193388 & 13.1111126 & 17.087418 & 13.048634 & 0.5320131 & 0\\
3 & 39.9510329 & 13.0036249 & 17.2513067 & 12.9993049 & 0.2760227 & 0\\
4 & 40.1496371 & 12.9966526 & 17.307515 & 12.9960858 & 0.0799361 & 0\\
5 & 40.2052568 & 12.9972241 & 17.3246431 & 12.9970895 & 0.0212722 & 0\\
6 & 40.2194854 & 12.9962604 & 17.3277278 & 12.996271 & 0.005506 & 0\\
7 & 40.2215552 & 12.9957037 & 17.3283368 & 12.9957233 & 0.0013865 & 0\\
8 & 40.2199042 & 12.995122 & 17.3271957 & 12.9951272 & 0.0003476 & 0\\
9 & 40.2189412 & 12.9970594 & 17.3293785 & 12.997061 & 8.73E-05 & 0\\
10 & 40.2212047 & 12.9959983 & 17.3285051 & 12.9959978 & 2.31E-05 & 0\\
11 & 40.2206441 & 12.9956756 & 17.3285992 & 12.9956754 & 6.00E-06 & 0\\
12 & 40.1675266 & 12.9963116 & 17.3284901 & 12.9963118 & 7.00E-07 & 0\\
13 & 40.7464532 & 12.9970054 & 17.3285922 & 12.9970032 & 7.00E-07 & 0\\
14 & 44.5853178 & 12.9967234 & 17.3279259 & 12.9967244 & 0 & 0\\
15 & 45.3009456 & 12.997981 & 17.3296289 & 12.99798 & 0 & 0\\
16 & 45.3487089 & 12.9976653 & 17.3282103 & 12.9976649 & 0 & 0\\
17 & 45.3496228 & 12.9961767 & 17.3260454 & 12.9961782 & 0 & 0\\
18 & 45.351359 & 12.9986208 & 17.3302856 & 12.9986212 & 0 & 0\\
19 & 45.3513136 & 12.9964465 & 17.3285827 & 12.9964469 & 0 & 0\\
20 & 45.3510285 & 12.9965624 & 17.3302195 & 12.9965623 & 0 & 0\\
21 & 45.3503113 & 12.9941995 & 17.3281136 & 12.9942002 & 0 & 0\\
22 & 45.3454522 & 12.9965376 & 17.3279103 & 12.996537 & 0 & 0\\
23 & 45.3151619 & 12.9971009 & 17.3284377 & 12.9971044 & 0 & 0\\
24 & 45.1841774 & 13.0202852 & 17.3391219 & 13.0202826 & 0 & 0\\
25 & 46.6619492 & 13.2828222 & 17.261221 & 13.2828215 & 0 & 0\\
26 & 58.5652028 & 16.8179593 & 21.0366657 & 13.7850471 & 0 & 0\\
27 & 0 & 4.0541284 & 4.3904588 & 2.1950835 & 0 & 0\\
28 & 0 & 0 & 0 & 0 & 0 & 0\\
29 & 0 & 0 & 0 & 0 & 0 & 0\\
    \end{tabular}
  \end{center}
\end{table}

\clearpage

Julia code for all of the algorithms proposed in this paper appear in Appendix A.

\section{Conclusions}

Both the naive (unfused) and naive (fused) algorithms are clear improvements on the {\tt Julia1.1\_hypot()} code. Moreover, with its more accurate calculation of $a^2+b^2$ the naive (fused) algorithm gives results that are directly comparable to the results from {\tt clib\_hypot()} with far less work. In many ways, this demonstrates the limits of trying to do better by just computing $a^2+b^2$ more accurately. However, the hands down winner in all of this is the corrected code (both versions). It requires substantially less work than the {\tt clib\_hypot()} code but delivers vastly superior results. 

\section*{Acknowledgments}
The author would like to thank Prof. Lucas Wilcox of the Naval Postgraduate School for invaluable help with the subtleties of the Julia programming language, and also Claude-Pierre Jeannerod of INRIA for several very helpful comments on the analysis.

\section*{Appendix A:  Julia Source Code for the Tested Algorithms}

\begin{verbatim}

function MyHypot1(x::T,y::T) where T<:AbstractFloat  # Naive (Unfused)
    #Return Inf if either or both imputs is Inf (Compliance with IEEE754)
    if isinf(x) || isinf(y)
        return convert(T,Inf)
    end

    # Order the operands
    ax,ay = abs(x), abs(y)
    if ay > ax
        ax,ay = ay,ax
    end

    # Widely varying operands
    if ay < ax*sqrt(eps(T)/2)  #Note: This also gets ay == 0
        return ax
    end

    # Operands do not vary widely
    scale = eps(sqrt(floatmin(T)))  #Rescaling constant
    if ax > sqrt(floatmax(T)/2)
        ax = ax*scale
        ay = ay*scale
        scale = inv(scale)
    elseif ay < sqrt(floatmin(T))
        ax = ax/scale
        ay = ay/scale
    else
        scale = one(scale)
    end
    sqrt(ax*ax+ay*ay)*scale
end

function MyHypot2(x::T,y::T) where T<:AbstractFloat # Naive (Fused)
    #Return Inf if either or both imputs is Inf (Compliance with IEEE754)
    if isinf(x) || isinf(y)
        return convert(T,Inf)
    end

    # Order the operands
    ax,ay = abs(x), abs(y)
    if ay > ax
        ax,ay = ay,ax
    end

    # Widely varying operands
    if ay <= ax*sqrt(eps(T)/2)  #Note: This also gets ay == 0
        return ax
    end

    # Operands do not vary widely
    scale = eps(sqrt(floatmin(T)))  #Rescaling constant
    if ax > sqrt(floatmax(T)/2)
        ax = ax*scale
        ay = ay*scale
        scale = inv(scale)
    elseif ay < sqrt(floatmin(T))
        ax = ax/scale
        ay = ay/scale
    else
        scale = one(scale)
    end
    sqrt(muladd(ax,ax,ay*ay))*scale
end

function MyHypot3(x::T,y::T) where T<:AbstractFloat  # Corrected (Unfused)
    #Return Inf if either or both imputs is Inf (Compliance with IEEE754)
    if isinf(x) || isinf(y)
        return convert(T,Inf)
    end

    # Order the operands
    ax,ay = abs(x), abs(y)
    if ay > ax
        ax,ay = ay,ax
    end

    # Widely varying operands
    if ay <= ax*sqrt(eps(T)/2)  #Note: This also gets ay == 0
        return ax
    end

    # Operands do not vary widely
    scale = eps(sqrt(floatmin(T)))  #Rescaling constant
    if ax > sqrt(floatmax(T)/2)
        ax = ax*scale
        ay = ay*scale
        scale = inv(scale)
    elseif ay < sqrt(floatmin(T))
        ax = ax/scale
        ay = ay/scale
    else
        scale = one(scale)
    end
    h = sqrt(ax*ax+ay*ay)

    # This is a well balanced single branch code
    # delta = h-ax
    #h -= ( delta*(2*(ax-ay)) + (2*delta-ay)*ay + delta*delta)/(2*h)
    # End single branch

    # This is a well balanced twopar branch code 
    if h <= 2*ay
        delta = h-ay
        h -= (ax*(2*delta - ax)+(delta-2*(ax-ay))*delta)/(2*h)
    else
        delta = h-ax
        h -= ( 2*delta*(ax-2*ay) + (4*delta-ay)*ay + delta*delta)/(2*h)
    end
    # End two branch code

    h*scale
end

function MyHypot4(x::T,y::T) where T<:AbstractFloat  # Corrected (Fused)
    #Return Inf if either or both imputs is Inf (Compliance with IEEE754)
    if isinf(x) || isinf(y)
        return convert(T,Inf)
    end

    # Order the operands
    ax,ay = abs(x), abs(y)
    if ay > ax
        ax,ay = ay,ax
    end

    # Widely varying operands
    if ay <= ax*sqrt(eps(T)/2)  #Note: This also gets ay == 0
        return ax
    end

    # Operands do not vary widely
    scale = eps(sqrt(floatmin(T)))  #Rescaling constant
    if ax > sqrt(floatmax(T)/2)
        ax = ax*scale
        ay = ay*scale
        scale = inv(scale)
    elseif ay < sqrt(floatmin(T))
        ax = ax/scale
        ay = ay/scale
    else
        scale = one(scale)
    end
    h = sqrt(fma(ax,ax,ay*ay))
    h_sq = h*h
    ax_sq = ax*ax
    x = fma(-ay,ay,h_sq-ax_sq) + fma(h,h,-h_sq) - fma(ax,ax,-ax_sq)
    h-=x/(2*h)
    h*scale
end


\end{verbatim}

\section*{Appendix B:  C Source Code for {\tt clib\_hypot()} }

\begin{verbatim}

/* @(#)e_hypot.c 1.3 95/01/18 */
/*
 * ====================================================
 * Copyright (C) 1993 by Sun Microsystems, Inc. All rights reserved.
 *
 * Developed at SunSoft, a Sun Microsystems, Inc. business.
 * Permission to use, copy, modify, and distribute this
 * software is freely granted, provided that this notice 
 * is preserved.
 * ====================================================
 */

#include "cdefs-compat.h"
//__FBSDID("$FreeBSD: src/lib/msun/src/e_hypot.c,v 1.14 2011/10/15 07:00:28 das Exp $");

/* __ieee754_hypot(x,y)
 *
 * Method :                  
 *	If (assume round-to-nearest) z=x*x+y*y 
 *	has error less than sqrt(2)/2 ulp, than 
 *	sqrt(z) has error less than 1 ulp (exercise).
 *
 *	So, compute sqrt(x*x+y*y) with some care as 
 *	follows to get the error below 1 ulp:
 *
 *	Assume x>y>0;
 *	(if possible, set rounding to round-to-nearest)
 *	1. if x > 2y  use
 *		x1*x1+(y*y+(x2*(x+x1))) for x*x+y*y
 *	where x1 = x with lower 32 bits cleared, x2 = x-x1; else
 *	2. if x <= 2y use
 *		t1*y1+((x-y)*(x-y)+(t1*y2+t2*y))
 *	where t1 = 2x with lower 32 bits cleared, t2 = 2x-t1, 
 *	y1= y with lower 32 bits chopped, y2 = y-y1.
 *		
 *	NOTE: scaling may be necessary if some argument is too 
 *	      large or too tiny
 *
 * Special cases:
 *	hypot(x,y) is INF if x or y is +INF or -INF; else
 *	hypot(x,y) is NAN if x or y is NAN.
 *
 * Accuracy:
 * 	hypot(x,y) returns sqrt(x^2+y^2) with error less 
 * 	than 1 ulps (units in the last place) 
 */

#include <float.h>
#include <openlibm_math.h>

#include "math_private.h"

OLM_DLLEXPORT double
__ieee754_hypot(double x, double y)
{
	double a,b,t1,t2,y1,y2,w;
	int32_t j,k,ha,hb;

	GET_HIGH_WORD(ha,x);
	ha &= 0x7fffffff;
	GET_HIGH_WORD(hb,y);
	hb &= 0x7fffffff;
	if(hb > ha) {a=y;b=x;j=ha; ha=hb;hb=j;} else {a=x;b=y;}
	a = fabs(a);
	b = fabs(b);
	if((ha-hb)>0x3c00000) {return a+b;} /* x/y > 2**60 */
	k=0;
	if(ha > 0x5f300000) {	/* a>2**500 */
	   if(ha >= 0x7ff00000) {	/* Inf or NaN */
	       u_int32_t low;
	       /* Use original arg order iff result is NaN; quieten sNaNs. */
	       w = fabs(x+0.0)-fabs(y+0.0);
	       GET_LOW_WORD(low,a);
	       if(((ha&0xfffff)|low)==0) w = a;
	       GET_LOW_WORD(low,b);
	       if(((hb^0x7ff00000)|low)==0) w = b;
	       return w;
	   }
	   /* scale a and b by 2**-600 */
	   ha -= 0x25800000; hb -= 0x25800000;	k += 600;
	   SET_HIGH_WORD(a,ha);
	   SET_HIGH_WORD(b,hb);
	}
	if(hb < 0x20b00000) {	/* b < 2**-500 */
	    if(hb <= 0x000fffff) {	/* subnormal b or 0 */
	        u_int32_t low;
		GET_LOW_WORD(low,b);
		if((hb|low)==0) return a;
		t1=0;
		SET_HIGH_WORD(t1,0x7fd00000);	/* t1=2^1022 */
		b *= t1;
		a *= t1;
		k -= 1022;
	    } else {		/* scale a and b by 2^600 */
	        ha += 0x25800000; 	/* a *= 2^600 */
		hb += 0x25800000;	/* b *= 2^600 */
		k -= 600;
		SET_HIGH_WORD(a,ha);
		SET_HIGH_WORD(b,hb);
	    }
	}
    /* medium size a and b */
	w = a-b;
	if (w>b) {
	    t1 = 0;
	    SET_HIGH_WORD(t1,ha);
	    t2 = a-t1;
	    w  = sqrt(t1*t1-(b*(-b)-t2*(a+t1)));
	} else {
	    a  = a+a;
	    y1 = 0;
	    SET_HIGH_WORD(y1,hb);
	    y2 = b - y1;
	    t1 = 0;
	    SET_HIGH_WORD(t1,ha+0x00100000);
	    t2 = a - t1;
	    w  = sqrt(t1*y1-(w*(-w)-(t1*y2+t2*b)));
	}
	if(k!=0) {
	    u_int32_t high;
	    t1 = 1.0;
	    GET_HIGH_WORD(high,t1);
	    SET_HIGH_WORD(t1,high+(k<<20));
	    return t1*w;
	} else return w;
}

#if LDBL_MANT_DIG == 53
__weak_reference(hypot, hypotl);
#endif
\end{verbatim}

\end{document}